\documentclass[11 pt]{article} 
\usepackage{amscd,amsmath}

\def\NN{\hbox{I\kern-.2em\hbox{N}}}
\def\TT{{\cdss\hbox{T\kern-.55em\hbox{T}}}}
\def\ZZ{{\cdss\hbox{Z\kern-.45em\hbox{Z}}}}
\def\RR{\hbox{I\kern-.2em\hbox{R}}}

\title{On sets represented by partitions}
\author{Jean-Christophe Aval}
\date{}

\font \petit=cmr9

\font \petitit=cmti9

\begin{document}

\maketitle

{\bf Abstract:} {\petit We prove a lemma that is useful to get upper bounds for the number of partitions without a given subsum. From this we can deduce an improved upper bound for the number of sets represented by the (unrestricted or into unequal parts) partitions of an integer {\petitit n}.}

\section{Introduction}

Let $n$ be an integer and let
$$n=n_1+n_2+...+n_j,\ n_i\in \NN^*,\ 1\le n_1\le n_2 \le ...\le n_j $$
be a partition $\Pi$ of $n$. We shall say that this partition represents an integer $a$ if there exist $\epsilon_1,\epsilon_2,...,\epsilon_j \in\{0,1\}$ such that $a=\sum_{i=1}^j \epsilon_i n_i$. Let ${\cal E}(\Pi)$ denote the set of these integers; we shall call it the set represented by $\Pi$.
One can easily see that ${\cal E}(\Pi)$ is included in $[0,n]$ and symmetric (if it contains $a$, it also contains $n-a$). 
For fixed $n$, let us introduce $p(n)$ the number of partitions of $n$ and $\hat p(n)$ the number of different sets amongst the sets ${\cal E}(\Pi)$ (where $\Pi$ runs over the $p(n)$ partitions of $n$).

Let $k$ be a positive integer. We shall say that a partition is $k$-reduced if and only if each summand appears at most $k$ times; for instance the 1-reduced partitions are the partitions into unequal parts.
We shall use $q(n,k)$ the number of $k$-reduced partitions of $n$ and $\hat q(n,k)$ the number of different ${\cal E}(\Pi)$ where $\Pi$ runs over the $q(n,k)$ $k$-reduced partitions of $n$. When $k$ equals 1, we shall note: $q(n)=q(n,1)$ and $\hat q(n)=\hat q(n,1)$. 

Following an idea due to P. Erd\"os, the sets represented by the partitions of an integer $n$ were first studied at the end of the 80's. P. Erd\"os, J.-L. Nicolas and A. S\'ark\"ozy (cf. [3]) obtained upper bounds for the number of partitions without a given subsum. P. Erd\"os then proposed to study the asymptotic behaviour of $\hat p(n)$ and $\hat q(n)$. In [1] and [5], M. Deléglise, P. Erd\"os, J.-L. Nicolas and A. S\'ark\"ozy proved the following estimates:

\vskip 0.2 cm
\bf
Theorem 1:
\it
For $n$ large enough, one has
$$q(n)^{0.51}\le\hat q(n)\le q(n)^{0.96}$$
and
$$p(n)^{0.361}\le\hat p(n)\le p(n)^{0.773}.$$
\rm
We shall obtain the following improved upper bounds:

\vskip 0.2 cm
\bf
Theorem 2:
\it
For $n$ large enough, one has
$$\hat q(n)\le q(n)^{0.955}\ \ {\it and}\ \ \hat p(n)\le p(n)^{0.768}.$$ 

\rm
To get these new exponents, we shall prove in part 2 a lemma improving a result due to J. Dixmier [2], whose application in part 3 gives the announced improvements.

\section{The main lemma}

Let $a$ be an integer, $a\le n$. We introduce ${\cal R}(n,a)$, the set of partitions of $n$ that do not represent $a$, and $R(n,a)$ shall denote its cardinality.
In the case of partitions into unequal parts, we shall need the same notions, with the similar notations ${\cal Q}(n,a)$ and $Q(n,a)$. We shall also define $Q(n,a,2)$ as the number of 2-reduced partitions $\Pi$ of $n$ such that $a$ is not represented by $\Pi$.

\bf
\vskip 0.2 cm
Lemma 1:
\it
Let $\epsilon>0.$
Assume there exists $\delta\in]0,1[$ such that, for any integer $n$ and for any integer $a$, the following property holds
$$\epsilon\sqrt n-1\le a\le 2\epsilon\sqrt n \Rightarrow R(n,a)\le p(n)^\delta.\leqno(1)$$
Then, for $n$ large enough, one has
$${{j}\over{2}}\epsilon \sqrt n\le a\le {{(j+1)}\over{2}}\epsilon \sqrt n \Rightarrow R(n,a)\le(2p(\epsilon\sqrt n))^{j-2}p(n)^\delta$$
\begin{itemize}
\item for $j=2,3,...,2[\sqrt n/2]$ if $\epsilon<1,$
\item for $j=2,3,...,\tau(n)$ with $\tau(n)=o(\sqrt n)$ for every $\epsilon$.
\end{itemize}

\vskip 0.2 cm
\rm
\bf Remark 1:
\rm
To obtain a similar conclusion, J. Dixmier [2] assumed that hypothesis (1) is true for $\epsilon\sqrt n\le a\le 3\epsilon\sqrt n$.

\vskip 0.2 cm
\bf
Proof:
\rm
We shall prove Lemma 1 by induction on $j$.
It is true for $j=2,3$ by (1). Let us suppose that $j\ge 4$ and that the result is true up to $j-1$. Let $a$ be such that ${{j}\over{2}}\epsilon \sqrt n\le a\le {{(j+1)}\over{2}}\epsilon \sqrt n$. Let $\Pi\in{\cal R}(n,a)$ and $b=[\epsilon\sqrt n]$. 

If $b$ is not represented by $\Pi$, then $\Pi$ belongs to a set ${\cal E}$ such that $|{\cal E}|\le p(n)^\delta.$

If $b$ is represented by $\Pi$, then we can write $\Pi=(\Pi',\Pi'')$, where $S(\Pi')=b$, $S(\Pi'')=n-b$ and $\Pi''$ does not represent $a-b$. We get
$$a-b\ge {{j}\over{2}}\epsilon\sqrt n-\epsilon\sqrt n={{j-2}\over{2}}\epsilon\sqrt n\ge \epsilon\sqrt {n-b}$$
since $j\ge 4$, and
$$a-b\le {{j+1}\over{2}}\epsilon\sqrt n-\epsilon\sqrt n+1={{j-1}\over{2}}\epsilon\sqrt n+1.$$
Moreover we have
$${{j}\over{2}}\epsilon\sqrt {n-b}\ge{{j}\over{2}}\epsilon\sqrt n \Bigl( 1-{{b}\over{n}} \Bigr)\ge{{j}\over{2}}\epsilon\sqrt n - \epsilon {{j}\over{2}}{{\epsilon\sqrt n}\over{\sqrt n}}.$$
We still have to show (at least for $n$ large enough)
$${{j-1}\over{2}}\epsilon\sqrt n+1\le{{j}\over{2}}\epsilon\sqrt n - {{j}\over{2}}\epsilon^2.$$
 
\begin{itemize}
\item If $\epsilon<1$, since $j/2\le \sqrt n/2$, we have to check the inequality 
$$-1/2\epsilon\sqrt n +1\le-\epsilon^2\sqrt n/2,$$
which is true when $n$ is large enough.

\item In the second case, we want to show 
$$-1/2\epsilon\sqrt n +1\le -1/2\epsilon^2{{\tau(n)}\over{\sqrt n}}.$$
This is true when $n$ is large enough by using the hypothesis on $\tau(n)$.

\end{itemize}

We finally get
$${{j}\over{2}}\epsilon\sqrt {n-b}\ge a-b.$$
We deduce from the induction hypothesis that $\Pi''$ belongs to a set ${\cal F}$ such that 
$$|{\cal F}|\le (2p(\epsilon\sqrt n))^{j-3}p(n)^\delta.$$
This implies that $\Pi$ belongs to a set ${\cal G}$ such that
$$|{\cal G}|\le p((\epsilon\sqrt n))(2p(\epsilon\sqrt n))^{j-3}p(n)^\delta.$$
Hence we have
$$R(n,a)\le p(n)^\delta + p((\epsilon\sqrt n))(2p(\epsilon\sqrt n))^{j-3}p(n)^\delta\le (2p(\epsilon\sqrt n))^{j-2}p(n)^\delta$$
which completes the proof of the lemma.

\vskip 0.2 cm
\bf
Remark 2:
\rm
It is easy to see that the result remains true when we replace all the $R(n,a)$'s by $Q(n,a)$'s or by $Q(n,a,2)$'s, i.e. when we deal with partitions into unequal parts or with 2-reduced partitions (in the proof, if $\Pi$ is into unequal parts, then $\Pi'$ and $\Pi''$ are also into unequal parts; the same phenomenon occurs when we are dealing with 2-reduced partitions).

\section{Applications} 

This lemma is useful to get upper bounds for $\hat p(n)$ and $\hat q(n)$ improving those obtained in [1]. Lemma 1 allows us to prove the following lemma: 

\vskip 0.2 cm
\bf
Lemma 2:
\it 
When $n\to \infty$ we have:
\begin{enumerate}
\item for $1.07\sqrt n\le a \le n-1.07\sqrt n$,
$$Q(n,a)\le\exp((1+o(1))1.732\sqrt n),$$
\item for $0.81\sqrt n\le a \le n-0.81\sqrt n$,
$$Q(n,a,2)\le\exp((1+o(1))1.969\sqrt n).$$
\end{enumerate}

\rm
To get Lemma 2 (the method is developped in [1]), we find upper bounds for $Q(n,a)$ and $Q(n,a,2)$ when $a$ ranges over the interval $[\epsilon\sqrt n,2\epsilon\sqrt n]$ and we choose the best $\epsilon$; then we use Lemma 1 and the results in [3].
  
From Lemma 2, we get Theorem 2 as in [1]. For instance, when studying $\hat q(n)$, we distinguish two cases according to whether the partition represents all integers between $1.07\sqrt n$ and $n-1.07 \sqrt n$ or not. We get this way
$$\hat q(n)\le n\exp((1+o(1))1.732\sqrt n)+2^{1.07\sqrt n}\le q(n)^{0.955}$$
since $q(n)=\exp((1+o(1))\pi\sqrt{n/3})$ (cf. [4]).

The method is the same for $\hat p(n)$, since $\hat p(n)=\hat q(n,2)$ [1, Théorème 1].

\vskip 0.2 cm
\bf
Remark 3:
\rm
The improvement on the exponents in the Theorem 2 is small ($5.10^{-3}$). This comes from the fact that the functions (cf. [1]) we bound on an interval $[x,2x]$ (and not $[x,3x]$, see Remark 1) have slow variations around their minimum value. Indeed, even replacing $[x,2x]$ by $[x,(1+\eta)x]$ with $\eta$ decreasing to $0$ would only lead to another small improvement ($4.10^{-3}$ less than our results). To make the exponents in the upper bounds really smaller, we need to find another method.

\rm
\section{References}

\begin{enumerate}

%\noindent
\item M. Deléglise, P. Erd\"os and J.-L. Nicolas, Sur les ensembles représentés
par les partitions d'un entier $n$, {\it Discr. Math.}, to appear.    
\vskip 0.2 cm

%\noindent
\item J. Dixmier, Partitions avec sous-sommes interdites, {\it Bull. Soc. Math.
Belgique}, {\bf 42} (1990), 477-500.
\vskip 0.2 cm

%\noindent
\item P. Erd\"os, J.-L. Nicolas and A. S\'ark\"ozy, On the number of partitions of $n$
without a given subsum II, {\it Analytic Number Theory}, edited by B. Berndt, H. Diamond, H. Halberstam, A. Hildebrand, Birkhaüser 1990, 207-236.
\vskip 0.2 cm

%\noindent
\item G. H. Hardy and S. Ramanujan, Asymptotic formulae in combinatory 
analysis, {\it Proc. London Math. Soc. (2)}, {\bf 17} (1918), 75-115.
\vskip 0.2 cm

%\noindent
\item J.-L. Nicolas and A. S\'ark\"ozy, On two partition problems, {\it Acta Math. 
Hung.}, {\bf 77} (1997), 95-121.

\end{enumerate}

\end{document}